\documentclass[twoside,12pt]{article}
\usepackage[]{amsmath,amsthm,amssymb}
\usepackage[latin1]{inputenc}

\begin{document}

\title{{\Large\bf  Discrete Kontorovich-Lebedev transforms}}

\author{Semyon  YAKUBOVICH}
\maketitle

\markboth{\rm \centerline{ Semyon   YAKUBOVICH}}{}
\markright{\rm \centerline{KONTOROVICH-LEBEDEV TRANSFORMS}}

\begin{abstract} {\noindent Discrete analogs of the classical Kontorovich-Lebedev transforms are introduced and investigated. It involves series with the modified Bessel function or Macdonald function $K_{in}(x), x >0, n \in \mathbb{N}, i $ is the imaginary unit, and incomplete Bessel functions. Several expansions of suitable  functions and sequences in terms of these series and integrals are established. As an application,  a Dirichlet boundary value problem in the upper half-plane  for inhomogeneous Helmholtz  equation is solved. }

\end{abstract}
\vspace{4mm}



\vspace{4mm}

\section {Introduction and preliminary results}

In 1946 N.N. Lebedev    proved  (cf. [2]) the  formula

$$f(\tau) = {2\over \pi^2} \ \tau \sinh(\pi\tau) \int_0^\infty {K_{i\tau} (x) \over x} \int_0^\infty K_{iy}(x) f(y) dy dx,\ \tau >0,\eqno(1.1)$$
which defines the Kontorovich-Lebedev transform in terms of the modified Bessel function or Macdonald function $K_\nu(z)$  given, for instance,  by the integrals
(see [3], Vol. I,  Entries 2.4.18.4,  2.5.54.6) 

$$K_{\nu}(z)= \int_{0}^\infty e^{-z\cosh (u) } \cosh(\nu u) du,\   {\rm Re}  z > 0, \  \nu \in \mathbb{C},\eqno(1.2)$$

$$K_{\nu}(x)= {1\over \cos(\pi \nu/2)} \int_0^\infty  \cos(x \sinh (u))  \cosh (\nu u) du,\quad x >0,\   -1 < {\rm Re} \nu < 1,\eqno(1.3)$$

$$K_{\nu}(x)= {1\over \sin(\pi \nu/2)} \int_0^\infty  \sin (x \sinh (u))  \sinh (\nu u) du,\quad x >0,\   -1 < {\rm Re} \nu < 1.\eqno(1.4)$$
 As  is known, the modified Bessel function $K_\nu(z)$ satisfies the ordinary differential equation
$$  z^2{d^2u\over dz^2}  + z{du\over dz} - (z^2+\nu^2)u = 0,\eqno(1.5)$$
for which it is the solution that remains bounded as $z$ tends to infinity on the real line. It has the asymptotic behavior [4]
$$ K_\nu(z) = \left( \frac{\pi}{2z} \right)^{1/2} e^{-z} [1+ O(1/z)], \qquad z \to \infty,\eqno(1.6)$$
and near the origin
$$ K_\nu(z) = O\left ( z^{-|{\rm Re}\nu|}\right), \ z \to 0,\eqno(1.7)$$
$$K_0(z) = -\log z + O(1), \ z \to 0. \eqno(1.8)$$
In the sequel we will use the Lebedev inequality for the modified Bessel function (see  [4], p.219)

$$ \left| K_{i\tau} (x) \right| \le A \ {x^{-1/4} \over \sqrt{ \sinh (\pi\tau)} },\quad x, \tau > 0,\eqno(1.9)$$
where $A > 0$ is an absolute constant.  

The Kontorovich-Lebedev transform is the basic index transform [4], involving integration with respect to the index (a parameter) of the modified Bessel function. 
Starting from the classical theory of Fourier series and Fourier transform, the problem of finding of discrete analogs for numerous integral transforms with special functions as kernels is quite important and interesting task in Analysis to create new methods for solutions to boundary value problems for PDE's.  In this paper we will undertake an attempt to introduce a discrete analog for the Kontorovich-Lebedev transform (1.1) and its various modifications.  We will see how it rises naturally, employing the theory of the incomplete Bessel functions [1] and Fourier series.  Our goals will involve new properties of these functions, their integral representations in order to establish several expansions of suitable classes of sequences and functions. As an application, a Dirichlet boundary value problem in the upper half-plane  for the inhomogeneous Helmholtz  equation will be solved explicitly.

\section{Incomplete Bessel functions, their properties and representations} 

Following [1],  we define the incomplete modified Bessel function $J(z,\nu,w)$, cutting  integral (1.2)

$$J(z,\nu,w) = \int_{0}^w  e^{-z\cosh (u) } \cosh(\nu u) du.\eqno(2.1)$$
On the other hand, Entry 2.4.18.13 in [3], Vol. I suggests the representation of  Macdonald's  function  in the form

$$K_\nu(x)= {1\over 2}\   e^{\nu\pi i/2} \int_{-\infty}^\infty e^{ix \sinh(u) -\nu u} du,\quad x >0,\  |{\rm Re}  \nu | < 1.\eqno(2.2)$$
Letting $\nu= i\tau,\ \tau \in \mathbb{R}$ in (2.2), we easily find the formula

$$K_{i\tau}(x)=  e^{-\pi \tau/2} \int_{0}^\infty \cos\left(x \sinh(u) - \tau u\right) du,\eqno(2.3)$$
which can be also verified via (1.3), (1.4).  In the sequel we will employ the incomplete versions of  formulas (1.3), (1.4), introducing the following functions

$$K_c\left(x,\nu,w\right)= {1\over \cos(\pi \nu/2)} \int_0^w  \cos(x \sinh (u))  \cosh (\nu u) du,\eqno(2.4)$$

$$K_s\left(x,\nu,w\right)= {1\over \sin(\pi \nu/2)} \int_0^w  \sin (x \sinh (u))  \sinh (\nu u) du.\eqno(2.5)$$
Further, the incomplete modified Bessel function (2.1) is a solution of the inhomogeneous differential equation (cf. [1])

$$  z^2{d^2u\over dz^2}  + z{du\over dz} - (z^2+\nu^2)u = - \left(\nu \sinh(\nu w) + z \cosh(\nu w)\sinh(w) \right) e^{-z\cosh(w)}.\eqno(2.6)$$

In particular,  letting  $w=\pi, \ \nu = in,\ n \in \mathbb{N}$ in (2.1), (2.4) and integrating  by parts, we obtain, respectively,

$$ J(x, in, \pi) =  {x\over n} \int_0^\pi e^{-x\cosh(u)} \sinh(u) \sin (n u) du,\quad x  >0,\eqno(2.7)$$

$$K_c (x, in, \pi)= {x\over n \cosh(\pi n/2)} \int_0^\pi \sin (x \sinh (u))  \cosh(u) \sin (n u) du,\quad   x   >0.\eqno(2.8)$$
Moreover, the theory of Fourier series suggests  the following equalities for all $x >0$ and $u \in [0,\pi]$

$${\pi\over 2} e^{-x\cosh(u)} =  {1\over 2}\ J (x,0, \pi)   + \sum_{n=1}^\infty  J(x, in, \pi) \cos(nu),\eqno(2.9)$$

$${\pi x\over 2 } e^{-x\cosh(u)} \sinh(u) =  \sum_{n=1}^\infty  n  J (x, in, \pi) \sin(nu),\eqno(2.10)$$

$${\pi\over 2} \cos(x\sinh(u))  =  {1\over 2}\  K_c(x, 0, \pi)  + \sum_{n=1}^\infty  \cosh\left({\pi n\over 2}\right)  K_c(x, in, \pi) \cos(nu),\eqno(2.11)$$

$$ {\pi x\over 2 } \sin (x \sinh (u))  \cosh(u) = \sum_{n=1}^\infty  n  \cosh\left({\pi n\over 2}\right) K_c(x, in, \pi)  \sin(nu),\eqno(2.12)$$

$${\pi\over 2} \sin(x\sinh(u))  =   \sum_{n=1}^\infty   \sinh\left({\pi n\over 2}\right)  K_s(x, in, \pi)  \sin(nu),\eqno(2.13)$$

$${\pi\over 2}\int_0^\pi  e^{-2 x\cosh(u)} du  =   {\pi\over 2}  \ J(2x, 0, \pi) =  {1\over 2}\  J^2 (x,0,\pi)  + \sum_{n=1}^\infty  J^2 (x, in, \pi),\eqno(2.14)$$

$${\pi x^2\over 2 }\int_0^\pi  e^{-2 x\cosh(u)} \sinh^2 (u) du  =   {\pi x^2 \over 4} \left[ J (2x, 2, \pi) -  J (2x,0, \pi) \right] $$

$$= \sum_{n=1}^\infty   n^2 J^2(x, in, \pi),\eqno(2.15)$$

$${\pi\over 2}  \int_0^\pi \cos^2 (x\sinh(u)) du   =  {\pi\over 4}  \left[ \pi +  K_c (2x, 0, \pi) \right] $$

$$= {\left[ K_c (x, 0, \pi)\right]^2  \over 2}  + \sum_{n=1}^\infty  \cosh^2\left({\pi n\over 2}\right)  \left[ K_c (x,in, \pi) \right]^2 ,\eqno(2.16)$$

$$ {\pi x^2\over 2 }  \int_0^\pi \sin^2 (x \sinh (u))  \cosh^2(u) du  = \sum_{n=1}^\infty  n^2  \cosh^2\left({\pi n\over 2}\right) \left[ K_c(x, in, \pi) \right]^2,\eqno(2.17)$$

$$ {\pi \over 2 }  \int_0^\pi \sin^2 (x \sinh (u)) du  = \sum_{n=1}^\infty    \sinh^2\left({\pi n\over 2}\right) \left[ K_s(x, in, \pi) \right]^2.\eqno(2.18)$$

The following lemma establishes biorthogonality of sequences of functions 

$$ \{ K_{in}(x) \}_{ n\in \mathbb{N}},  \{ J (x, in, \pi) \}_{ n\in \mathbb{N}},\eqno(2.19)$$

$$\{ K_{in}(x) \}_{ n\in \mathbb{N}}, \  \{ K_c (x, in, \pi) \}_{ n\in \mathbb{N}},\eqno(2.20) $$

$$ \{ K_{i\tau}(n) \}_{ n\in \mathbb{N}},\  \{ K_s \left(n, i\tau, \sinh^{-1} (\pi)\right) \}_{ n\in \mathbb{N}},\eqno(2.21)$$
where we use the notation $ \sinh^{-1}(u) = \log\left( u+ \sqrt{1+ u^2}\right).$

Indeed, we have

{\bf Lemma 1}. {\it  Sequences of functions $(2.19), (2.20)$  are biorthogonal over $\mathbb{R}_+$ with respect to the measure ${dx\over x}$, i.e.

$$I_1= \int_0^\infty  K_{in}(x) J(x,im, \pi) {dx\over x} = { \pi^2\over 2} \ {  \delta_{n,m}\over n \sinh(\pi n)},\quad  n, m \in \mathbb{N},\eqno(2.22)$$

$$I_2= \int_0^\infty  K_{in}(x) K_c(x, im, \pi) {dx\over x} = { \pi^2\over 2} \ {  \delta_{n,m}\over n \sinh(\pi n) },\quad  n, m \in \mathbb{N},\eqno(2.23)$$
where $\delta_{n,m}$ is the Kronecker symbol.  Besides, sequences $(2.21)$ are biorthogonal over $\mathbb{R}_+$ with respect to the measure $\tau \sinh(\pi\tau) d\tau$, namely,}

$$I_3= \int_0^\infty   K_{i\tau}(n)  K_s\left(m, i\tau, \sinh^{-1} (\pi)\right) \tau \sinh(\pi\tau) d\tau =  { \pi^2 n \over 2} \  \delta_{n,m},\quad  n, m \in \mathbb{N}.\eqno(2.24)$$

\begin{proof}   In fact, using (2.7), we substitute the corresponding integral into (2.22), having the estimate

$$ |I_1| \le {1\over m} \int_0^\infty  \left| K_{in}(x) \right| \int_0^\pi e^{-x\cosh(u)} \sinh(u)  \left| \sin (m u) \right| du dx $$

$$\le  {1\over m} \int_0^\infty  K_{0}(x)  \int_0^\pi e^{-x\cosh(u)} \sinh(u)  du dx =  {1\over m} \int_0^\infty  K_{0}(x) \left( 1- e^{-x\cosh(\pi)}\right) {dx \over x}  < \infty.$$
Therefore Fubini's theorem permits the interchange of the order of integration to write $I_1$  as follows

$$I_1 =   {1\over m} \int_0^\pi  \sinh(u)  \sin (m u)   \int_0^\infty  K_{in}(x) e^{-x\cosh(u)} dx du.\eqno(2.25)$$
But the   integral with respect to $x$ on the right-hand side of (2.25) is calculated in [3], Vol. II, Entry 2.16.6.1,  and we have

$$ \int_0^\infty  e^{-x\cosh(u)} K_{in}(x)  dx = {\pi \sin( nu) \over \sinh (u) \sinh(\pi n)}.\eqno(2.26)$$
Hence (2.25) becomes

$$I_1 =   {\pi \over m \sinh(\pi n)} \int_0^\pi  \sin (n u)   \sin( mu) du =  { \pi^2\over 2} \ {  \delta_{n,m}\over n \sinh(\pi n)},$$
which proves (2.22).   

Analogously, employing (2.8) and the estimate 

$$ |I_2| \le  {1\over m \cosh(\pi m/2)} \int_0^\infty  \left| K_{in}(x) \right|   \int_0^\pi \left| \sin (x \sinh (u))  \cosh(u) \sin (m u) \right| du  dx $$

$$\le   {\pi \cosh(\pi) \over m \cosh(\pi m/2)} \int_0^\infty   K_{0}(x) dx < \infty,$$
we change the order of integration by Fubini's theorem to obtain $I_2$  in the form 

$$I_2  =  {1\over m \cosh(\pi m/2)}  \int_0^\pi   \cosh(u) \sin (m u)   \int_0^\infty  K_{in}(x) \sin (x \sinh (u)) dx du.$$
The latter integral with respect to $x$ is calculated in   [3], Vol. II, Entry   2.16.14.1, and we get

$$ \int_0^\infty  K_{in}(x) \sin (x \sinh (u)) dx =  {\pi\over 2}  { \sin (n u)\over \sinh(\pi n/2) \cosh(u) }.\eqno(2.27)$$
Therefore

$$I_2=  {\pi \over 2 m \cosh(\pi m/2)  \sinh(\pi n/2)}  \int_0^\pi  \sin (n u)  \sin (m u) du =    {\pi^2 \delta_{n,m}  \over 2 n \sinh(\pi n)} $$
and we establish equality (2.23).  

Finally, in order to prove (2.24) we appeal to (2.5) to write

$$K_s \left(m, i\tau, \sinh^{-1} (\pi) \right)= {1\over \sinh(\pi \tau/2)} \int_0^{\sinh^{-1} (\pi) } \sin(m \sinh (u))  \sin (\tau u) du.$$
Integration by parts and simple substitutions on the right-hand side of the latter equality then give

$${1\over \sinh(\pi \tau/2)} \int_0^{\sinh^{-1} (\pi) } \sin(m \sinh (u))  \sin (\tau u) du =  {m \over \tau \sinh(\pi \tau/2)} \int_0^{\sinh^{-1} (\pi) } \cos(m \sinh (u)) $$

$$\times  \cos (\tau u) \cosh(u) du =  {m \over \tau \sinh(\pi \tau/2)} \int_0^{\pi} \cos(m u)  \cos \left(\tau \sinh^{-1} (u) \right) du.$$
Hence

$$I_3 = 2m \int_0^\infty   K_{i\tau}(n)  \cosh \left({\pi\tau\over 2}\right)  \int_0^{\pi} \cos(m u)  \cos \left(\tau \sinh^{-1} (u)\right) du d\tau$$

$$= 2m \lim_{N\to \infty } \int_0^N   K_{i\tau}(n)  \cosh \left({\pi\tau\over 2}\right)  \int_0^{\pi} \cos(m u)  \cos \left(\tau \sinh^{-1} (u)\right) du d\tau$$

$$= 2m \lim_{N\to \infty }   \int_0^{\pi} \cos(m u)  \int_0^N   K_{i\tau}(n)  \cosh \left({\pi\tau\over 2}\right)  \cos \left(\tau \sinh^{-1} (u)\right) d\tau du,\eqno(2.28)$$
where the latter interchange of the order of integration is due to the absolute and uniform convergence.  The problem now is to motivate the passage to the limit under the integral sign on the right-hand side of the latter equality. Then we use the value of the improper integral with respect to $\tau$ (cf. [4]) 

 $$ \int_0^\infty   K_{i\tau}(n)  \cosh \left({\pi\tau\over 2}\right)  \cos\left(\tau  \sinh^{-1} (u) \right) d\tau = {\pi\over 2} \cos(n u)\eqno(2.29)$$  
to complete the proof of Lemma 1.   So, choosing a big enough $N >0$ and fixing $n \in \mathbb{N}$,  we use asymptotic behavior of the Macdonald function with respect to $\tau$ [4] 

$$ K_{i\tau}(n) = \left({2\pi\over \tau}\right)^{1/2} e^{-\pi\tau/2} \sin\left( {\pi\over 4} + \tau \log\left({2\tau\over n} \right) - \tau \right) \left( 1+ O\left({1\over \tau}\right) \right),\ \tau \to +\infty$$
to get the remainder of the integral (2.28) in the form

$$ \int_N^\infty   K_{i\tau}(n)  \cosh \left({\pi\tau\over 2}\right)  \cos\left(\tau  \sinh^{-1} (u) \right) d\tau $$

$$= \sqrt{{\pi\over 2}}  \int_N^\infty   \sin\left( {\pi\over 4} + \tau \log\left({2\tau\over n} \right) - \tau \right)  \cos\left(\tau \sinh^{-1} (u)\right) {d\tau\over \sqrt \tau }
+  O\left({1\over \sqrt N}\right). $$
As we see via trigonometric identities the latter integral is associated with the following ones

$$G_1(N)= \int_N^\infty   \sin\left(  \tau \left( \log\left({2\tau\over n} \right) - 1 \mp \sinh^{-1} (u) \right)\right) {d\tau\over \sqrt \tau }, $$

$$G_2(N)= \int_N^\infty   \cos\left(  \tau \left( \log\left({2\tau\over n} \right) - 1 \mp \sinh^{-1} (u) \right) \right) {d\tau\over \sqrt \tau }.$$
Taking, for instance, the sine integral, we integrate by parts to derive

$$ G_1(N)=  {1\over \sqrt N} \ \frac{ \cos\left(  N \left( \log\left({2N/ n} \right) - 1 \mp \sinh^{-1} (u) \right)\right)}{ \log( 2N /n) \mp \sinh^{-1} (u) }$$

$$- {1\over 2} \int_N^\infty  { \cos \left(  \tau \left( \log\left({2\tau/ n} \right) - 1 \mp \sinh^{-1} (u) \right)\right) \over \left(\log(2\tau/n) \mp \sinh^{-1} (u) \right) \tau^{3/2}  } d\tau$$

$$-  \int_N^\infty  { \cos \left(  \tau \left( \log\left({2\tau/ n} \right) - 1 \mp \sinh^{-1} (u) \right)\right) \over \left(\log(2\tau/n) \mp \sinh^{-1} (u)  \right)^2 \tau^{3/2}  } d\tau.$$
Hence,

$$| G_1(N)| \le   { 1 \over  \sqrt N}\  \left( \log( 2N /n) - \sinh^{-1} (\pi)  \right)^{-1}  $$

$$ +   \int_N^\infty {\left(\log(2\tau/n) -  \sinh^{-1} (\pi) \right) /2 +1\over  \left(\log(2\tau/n) - \sinh^{-1} (\pi)  \right)^2 \tau^{3/2}  } d\tau \to 0,\ N \to \infty.$$
Analogously the integral $G_2(N)$ can be treated.  Thus the passage to the limit under the integral sign in (2.28) is allowed. Then in view of formula (2.29), we readily end up with the biorthogonality (2.24).    

\end{proof} 

To complete this section, we will exhibit below our main results, which will be proved in the sequel. It involves  the following expansions  of   sequences and functions under certain conditions,  giving  rise to  discrete analogs of the Kontorovich - Lebedev transform (1.1), namely
$( n\in  \mathbb{N},\ x >0 )$

$$a_n =  {2\over \pi^2} \ n  \sinh(\pi n) \int_0^\infty J (x, in, \pi) \sum_{m=1}^\infty a_m K_{im}(x)  {dx\over x} ,\eqno(2.30)$$

$$a_n =  {2\over \pi^2} \ n  \sinh(\pi n) \int_0^\infty K_{in} (x)  \sum_{m=1}^\infty a_m J (x, im, \pi)  {dx\over x},\eqno(2.31)$$

$$a_n =  {2\over \pi^2} \ n  \sinh(\pi n) \int_0^\infty K_c(x, in,  \pi)  \sum_{m=1}^\infty a_m K_{im}(x)  {dx\over x},\eqno(2.32)$$

$$a_n =  {2\over \pi^2} \ n  \sinh(\pi n) \int_0^\infty K_{in} (x)  \sum_{m=1}^\infty a_m K_c (x, im, \pi)  {dx\over x},\eqno(2.33)$$

$$a_n =  {4\over \pi^2\ n} \lim_{\alpha\to  {\pi\over 2} -} \int_0^\infty   \tau  \sinh\left({\pi \tau\over 2}\right) \cosh\left(\alpha \tau \right)  K_s\left(n, i\tau, \sinh^{-1} (\pi) \right) \sum_{m=1}^\infty a_m K_{i\tau}(m)  d\tau,\eqno(2.34)$$

$$a_n =  {2\over \pi^2\ n}  \int_0^\infty   \tau  \sinh(\pi \tau )    K_{i\tau}(n) \sum_{m=1}^\infty a_m K_s \left(m, i\tau, \sinh^{-1} (\pi)  \right)  d\tau,\eqno(2.35)$$

$$  f(x)  =  {2\over x \pi^2} \   \sum_{n=1}^\infty   n  \sinh(\pi n)  J (x, in, \pi)  \int_0^\infty  K_{in} (y) f(y) dy,\eqno(2.36)$$

$$f(x)  =  {2\over \pi^2} \   \sum_{n=1}^\infty   n  \sinh(\pi n)  K_{in} (x)   \int_0^\infty  J (y, in, \pi) f(y) {dy\over y},\eqno(2.37)$$

$$  f(x)  =  {2\over x \pi^2} \   \sum_{n=1}^\infty   n  \sinh(\pi n)  K_c(x, in,  \pi)  \int_0^\infty  K_{in} (y)  f(y) dy,\eqno(2.38)$$

$$f(x)  =  {2\over \pi^2} \   \sum_{n=1}^\infty   n  \sinh(\pi n)  K_{in} (x)  \int_0^\infty K_c (y, in, \pi) f(y) {dy\over y},\eqno(2.39)$$

$$f(x)  =  {4\over \pi^2\ } \  x \sinh\left({\pi x\over 2}\right)    \sum_{n=1}^\infty  {1\over n} \ K_s \left(n, ix,  \sinh^{-1}(\pi) \right) $$

$$\times \int_0^\infty   \cosh\left({\pi \tau\over 2}\right )  K_{i\tau}(n)  f(\tau) d\tau,\eqno(2.40)$$

$$f(x)  =  {2\over \pi^2\ }   \sum_{n=1}^\infty  { K_{ix}(n) \over n} \int_0^\infty   \tau \sinh\left(\pi \tau\right )  K_s\left(n, i\tau,  \sinh^{-1}(\pi)\right)   f(\tau) d\tau.\eqno(2.41)$$

\section {Expansion theorems}

We begin with

{\bf Theorem 1.} {\it Let the sequence $\{a_n\}_{n\in \mathbb{N}}$ be such that the following series converges

$$\sum_{n=1}^\infty  |a_n|  e^{-\pi n/2} < \infty.\eqno(3.1)$$
Then  expansion $(2.30)$ holds, where the iterated series and integral converge absolutely.}

\begin{proof}  In fact, recalling representation (2.7) of the kernel $J(x, in, \pi)$, substituting this expression on  the right-hand side of (2.30) and employing Lebedev's inequality (1.9),  we derive the estimates 

$$ {2\over \pi^2}  \sinh(\pi n) \int_0^\infty    \int_0^\pi \left| e^{-x\cosh(u)} \sinh(u) \sin (n u) \right| du  \sum_{m=1}^\infty \left| a_m K_{im}(x) \right| dx $$

$$\le  {2 A\over \pi^2}  \sinh(\pi n) \int_0^\infty   x^{-1/4}  \int_0^\pi e^{-x\cosh(u)} \sinh(u) du dx  \sum_{m=1}^\infty  {| a_m | \over \sqrt{ \sinh (\pi m)} } $$

$$  \le {2 \sqrt 2\over \pi^2} \  \Gamma(3/4) A \sinh(\pi n) \int_0^\pi { \sinh(u) \over \cosh^{3/4} (u)} du \sum_{m=1}^\infty  |a_m|  e^{-\pi m/2} $$

$$ =   {8 \sqrt 2\over \pi^2} \    \Gamma(3/4) A \sinh(\pi n) \left[ \cosh^{1/4} (\pi) -1\right]  \sum_{m=1}^\infty  |a_m|  e^{-\pi m/2} < \infty $$
owing to the condition (3.1) of the theorem. Hence the interchange of the order of integration and summation in (2.30) is guaranteed by Fubini's theorem and we find

$$ {2\over \pi^2}  \sinh(\pi n) \int_0^\infty   \int_0^\pi e^{-x\cosh(u)} \sinh(u) \sin (n u) du \sum_{m=1}^\infty  a_m K_{im}(x)  dx $$

$$ = {2\over \pi^2}  \sinh(\pi n) \sum_{m=1}^\infty  a_m   \int_0^\pi  \sinh(u) \sin (n u)  \int_0^\infty  e^{-x\cosh(u)} K_{im}(x)  dx du.$$
Then appealing to formula (2.26),  we have finally

$${2\over \pi^2}  \sinh(\pi n) \sum_{m=1}^\infty  a_m   \int_0^\pi  \sinh(u) \sin (n u)  \int_0^\infty  e^{-x\cosh(u)} K_{im}(x)  dx du$$

$$=  {2\over \pi}  \sinh(\pi n) \sum_{m=1}^\infty  {a_m \over \sinh(\pi m)}  \int_0^\pi  \sin (n u) \sin (m u) du = a_n.$$
Theorem 1 is proved. 

\end{proof}

Concerning expansion (2.31), we have 

{\bf Theorem 2.} {\it Let the sequence $\{a_n\}_{n\in \mathbb{N}}$ be such that the following series converges

$$\sum_{n=1}^\infty  {|a_n| \over n}  < \infty.\eqno(3.2)$$
Then  expansion $ (2.31)$ holds valid, where the iterated series and integral converge absolutely.}

\begin{proof}  Doing in the same manner, we appeal to (1.9) and (2.7)  to get the estimate

$$ {2\over \pi^2}  n \sinh(\pi n) \int_0^\infty  \left| K_{in}(x) \right|  \sum_{m=1}^\infty  { |a_m| \over m}  \int_0^\pi  \left| e^{-x\cosh(u)} \sinh(u) \sin (m u)  \right| du   dx $$

$$\le  {2 A\over \pi^2} n \sqrt{\sinh(\pi n)} \int_0^\infty   x^{-1/4}  \int_0^\pi e^{-x\cosh(u)} \sinh(u) du dx  \sum_{m=1}^\infty  {| a_m | \over m } $$

$$  =  {2  \Gamma(3/4) A\over \pi^2} n  \sinh(\pi n) \int_0^\pi { \sinh(u) \over \cosh^{3/4} (u)} du \sum_{m=1}^\infty { |a_m| \over m}  $$

$$ =   {8  \Gamma(3/4) A\over \pi^2}  n \sinh(\pi n) \left[ \cosh^{1/4} (\pi) -1\right]  \sum_{m=1}^\infty  {|a_m| \over m} < \infty $$
due to the condition (3.2) of the theorem. Therefore, as above,  we change  the order of integration and summation by Fubini's theorem and recall  equality (2.26) to  arrive at (2.31), completing the proof of Theorem 2.

\end{proof} 

{\bf Remark 1}.  Expansions (2.30), (2.31) generate two pairs of discrete Kontorovich-Lebedev transforms

$$f(x)= \sum_{m=1}^\infty  a_m K_{im} (x),\quad  x >0,\eqno(3.3)$$

 $$g(x)= \sum_{m=1}^\infty  b_m J (x,im, \pi),\quad  x >0.\eqno(3.4)$$
Coefficients $a_n,\ b_n,\ n \in \mathbb{N}$ are calculated accordingly

$$a_n =  {2\over \pi^2} \ n  \sinh(\pi n) \int_0^\infty J (x,in, \pi)  f(x)   {dx \over x},\eqno(3.5)$$

$$b_n =  {2\over \pi^2} \ n  \sinh(\pi n) \int_0^\infty K_{in} (x)   g(x)   {dx \over x}.\eqno(3.6)$$
Integral  (3.6) is the discrete  Kontorovich-Lebedev transform  whose existence conditions are well known in different spaces of functions (see [4]).  A simple conditions is based on the inequality $ | K_{in} (x)  | \le K_0 (x),\ x > 0$. Thus if $g$ is Lebesgue integrable with respect to the measure $K_0(x) dx/x$, i.e. $ g \in L_1\left(\mathbb{R}_+;    K_0(x) dx/x \right)$, the Kontorovich-Lebedev transform (3.6) exists.    The corresponding existence conditions for discrete transform (3.5) follow from the behavior of the kernel $J (x,in, \pi) $.  For instance, representation (2.7) yields the equality $J (x,in, \pi) = O(x)$ as $ x \to 0$ and  simple inequality

$$\left| J (x,in, \pi) \right| \le {x\over n} \int_0^\pi e^{-x\cosh(u)} \sinh(u) du =   {1\over n} \left( e^{-x}- e^{- x\cosh (\pi)}\right) \le  {2\over n}\   e^{-x},\ x > 0.$$ 
Hence,  if $f \in L_1\left(\mathbb{R}_+;    e^{-x} dx \right)$, the discrete transform (3.5) exists. 

{\bf Remark 2.}  Another way to prove the validity of expansions (2.30), (2.31) is to use the biorthogonality (2.22), verifying first   the interchange of the order of integration and summation under conditions (3.1), (3.2), respectively. 

{\bf Theorem 3}. {\it Let the sequence $\{a_n\}_{n\in \mathbb{N}}$ satisfy the condition

$$ \sum_{m=1}^\infty \left| a_m \right| e^{-\delta m} < \infty,\ \delta \in \left[0, \  {\pi\over 2} \right).\eqno(3.7)$$
Then expansions  $(2.32), (2.33)$ hold, where the iterated   series and integral converge absolutely.   Moreover,  expansion $(2.33)$  holds for a wider class of sequences under the condition }

$$ \sum_{m=1}^\infty {\left| a_m \right|\over m}  e^{-\pi m/2} < \infty.\eqno(3.8)$$

\begin{proof}  Indeed, appealing to the inequality for the modified Bessel function (see [4], p. 15)

$$ \left| K_{i\tau} (x) \right| \le e^{-\delta \tau} K_0\left( x \cos(\delta) \right), \quad x, \tau >0,\ \delta \in \  \left[0, \  {\pi\over 2} \right),$$
we use integral representation (2.8) for the kernel $K_c(x,in, \pi)$ and the following estimates to justify the interchange of the order of integration and summation in (2.32), ( 2.33) under conditions of the theorem.  Thus

$$\int_0^\infty {1\over x} \  | K_c (x,in, \pi) | \sum_{m=1}^\infty \left| a_m K_{im}(x) \right| dx \le  {1\over n \cosh(\pi n/2)}  \int_0^\infty  K_0\left( x \cos(\delta) \right) dx $$

$$\times  \int_0^\pi \left| \sin (x \sinh (u))  \cosh(u) \sin (n u)\right| du  \sum_{m=1}^\infty \left| a_m  \right|   e^{-\delta m}$$

$$\le  {\sinh(\pi) \over n \cosh(\pi n/2)}  \int_0^\infty  K_0\left( x \cos(\delta) \right) dx  \sum_{m=1}^\infty \left| a_m  \right|   e^{-\delta m} $$

$$= { \pi \sinh(\pi) \over 2 n \cosh(\pi n/2) \cos(\delta)}   \sum_{m=1}^\infty \left| a_m  \right|   e^{-\delta m}< \infty, \quad \delta \in \  \left[0, \  {\pi\over 2} \right),$$

$$\int_0^\infty { | K_{in}(x) |\over x} \sum_{m=1}^\infty \left| a_m K_c (x, im, \pi)  \right| dx $$

$$\le 2  \sinh(\pi) \int_0^\infty  K_0\left( x  \right) dx  \sum_{m=1}^\infty {\left| a_m \right| \over m }  e^{-\pi m/2}$$

$$\le \pi  \sinh(\pi) \sum_{m=1}^\infty \left| a_m \right|   e^{- \delta m}< \infty,\quad \delta \in \  \left[0, \  {\pi\over 2} \right).$$
Therefore, changing the order of integration and summation on the right-hand side of (2.32), (2.33) and employing biorthogonality  (2.23), we establish the desired equalities. 

\end{proof}

Concerning expansions (2.34), (2.35), we have 

{\bf Theorem 4}. {\it  Let the sequence $\{a_n\}_{n \in \mathbb{N} } \in l_1$, i.e. $\sum_{m\ge 1} |a_m| < \infty$. Then  its general term  can be expanded with respect to  $(2.33)$ with  the convergence  in the sense of Abel.  Moreover, under  stronger  condition
$$\sum_{m=1}^\infty  m \left| a_m \right|  < \infty\eqno(3.9)$$
expansion $(2.35)$  holds, where the convergence is in the improper sense.} 

\begin{proof} In fact, from the proof of Lemma 1 we have the equality 

$$ K_s\left(n, i\tau, \sinh^{-1}(\pi) \right) =   {n \over \tau \sinh(\pi \tau/2)}  \int_0^{\pi} \cos(n u)  \cos \left(\tau  \sinh^{-1}(u)\right) du.\eqno(3.10)$$
Substituting this expression on the right-hand side of (2.34), we get

$${4\over \pi^2\ n} \lim_{\alpha \to {\pi\over 2} -}  \int_0^\infty   \tau  \sinh\left({\pi \tau\over 2}\right )\cosh \left(\alpha \tau\right)   K_s\left(n, i\tau, \sinh^{-1}(\pi) \right)$$

$$\times  \sum_{m=1}^\infty a_m K_{i\tau}(m)  d\tau = {4\over \pi^2} \lim_{\alpha \to {\pi\over 2} -}  \int_0^\infty   \cosh \left(\alpha \tau\right)  \int_0^{\pi} \cos(n u)  \cos \left(\tau \sinh^{-1}(u)  \right) du $$

$$\times  \sum_{m=1}^\infty a_m K_{i\tau}(m)  d\tau.\eqno(3.11)$$
From  inequality (1.9) we have the estimate

$$\int_0^\infty   \cosh \left(\alpha \tau\right)  \int_0^{\pi} \left| \cos(n u)  \cos \left(\tau \sinh^{-1}(u) \right)\right| du   \sum_{m=1}^\infty \left| a_m K_{i\tau}(m) \right| d\tau $$

$$\le \pi A  \int_0^\infty   {\cosh \left(\alpha \tau\right)  \over \sqrt{\sinh(\pi\tau)}} \ d\tau \sum_{m=1}^\infty   \left| a_m \right| < \infty,\quad  0 < \alpha < {\pi\over 2},$$
which validates  the interchange of the order of integration and summation on the right-hand side of (3.11). Thus we obtain 

$${ 4\over \pi^2} \lim_{\alpha \to {\pi\over 2} -}  \int_0^\infty   \cosh \left(\alpha \tau\right)  \int_0^{\pi} \cos(n u)  \cos \left(\tau \sinh^{-1}(u) \right) du $$

$$\times  \sum_{m=1}^\infty a_m K_{i\tau}(m)  d\tau = { 4\over \pi^2} \lim_{\alpha \to {\pi\over 2} -}    \sum_{m=1}^\infty a_m  \int_0^{\pi} \cos(n u)  $$

$$\times \int_0^\infty   \cosh \left(\alpha \tau\right) K_{i\tau}(m)  \cos \left(\tau \sinh^{-1}(u) \right) d\tau du.$$
But the integral with respect to $\tau$ is calculated in [3], Vol. II, Entry 2.16.48.20 and we find the equality 

$$ \int_0^\infty   \cosh \left(\alpha \tau\right) K_{i\tau}(m)  \cos \left(\tau \sinh^{-1}(u) \right) d\tau = {\pi\over 2}\  e^{-m\cos(\alpha) (1+u^2)^{1/2}} \cos( m u\sin(\alpha)).\eqno(3.12)$$
Therefore, combining with (3.11), it gives 

$$ {4\over \pi^2\ n} \lim_{\alpha \to {\pi\over 2} -}  \int_0^\infty   \tau  \sinh\left({\pi \tau\over 2}\right )\cosh \left(\alpha \tau\right)   K_s\left(n, i\tau, \sinh^{-1}(\pi) \right)$$

$$\times  \sum_{m=1}^\infty a_m K_{i\tau}(m)  d\tau =  {2\over \pi} \lim_{\alpha \to {\pi\over 2} -} \sum_{m=1}^\infty a_m$$ 

$$\times   \int_0^{\pi}  e^{-m\cos(\alpha) (1+u^2)^{1/2}}  \cos(n u) \cos( m u\sin(\alpha)) du. \eqno(3.13)$$
Hence equality (2.34) will be proved,  if we motivate the passage to the limit under the series sign on the right-hand side of (3.13).  But this fact is an immediate consequence of the estimate

$$\sum_{m=1}^\infty \left| a_m \int_0^{\pi}  e^{-m\cos(\alpha) (1+u^2)^{1/2}}  \cos(n u) \cos( m u\sin(\alpha)) du \right| \le \pi \sum_{m=1}^\infty  \left| a_m \right|.$$
In order to establish  (2.35),  we appeal to the biorthogonality result (2.24) together with  (2.28) to justify changing   order of integration and summation.   The completion of the proof is left to the reader.

\end{proof}

Let us argue  the validity of expansions  (2.36), (2.37).   We have

{\bf Theorem 5}. {\it Let $f$ be a complex-valued function on $\mathbb{R}_+$ which is represented by the integral 

$$f(x)= \int_{-\pi}^\pi e^{-x\cosh(u)} \varphi(u) du,\ x >0,\eqno(3.14)$$
where $ \varphi(u) = \psi(u)\sinh(u)$ and $\psi$ is a  $2\pi$-periodic function, satisfying the Lipschitz condition on $[-\pi, \pi]$, i.e.

$$\left| \psi(u) - \psi(v)\right| \le C |u-v|, \quad  \forall \  u, v \in  [-\pi, \pi],\eqno(3.15)$$
where $C >0$ is an absolute constant.  Then expansion $(2.36)$ holds for all $x >0$.  Further, if $f$ can be expanded in terms of the series $(3.3)$
of a suitable   sequence $\{a_m\}_{m \in \mathbb{N}}$, satisfying condition $(3.1)$, then expansion $(2.37)$ holds valid. }

\begin{proof}  Let $S_N(x) $ denote  a partial sum of the series (2.36), substitute the value of $f(x)$ by integral (3.14) and  interchange the order of integration with the use of the equality (2.26) to obtain 

$$S_N(x)  = {2\over \pi  x} \   \sum_{n=1}^N   n \  J (x, in, \pi)  \int_{-\pi}^\pi   {\varphi(u) \over \sinh(u)} \sin(nu) du.\eqno(3.16)$$
But the right-hand side of (3.16) can be rewritten, employing representation (2.7) to get 

$$S_N(x)  =   {1\over \pi} \   \sum_{n=1}^N  \int_{-\pi}^\pi e^{-x\cosh(t)} \sinh(t) \sin (n t) dt  \int_{-\pi}^\pi   {\varphi(u) \over \sinh(u)} \sin(nu) du.\eqno(3.17)$$
Now, calculating the sum via the known formula, equality (3.17) becomes

$$  S_N(x)  =  {1\over 4 \pi} \   \int_{-\pi}^\pi e^{-x\cosh(t)} \sinh(t)   \int_{-\pi}^{\pi}   {\varphi(u)+ \varphi(-u)  \over \sinh(u)} \  {\sin \left((2N+1) (u-t)/2 \right)\over \sin( (u-t) /2)}  du dt $$

$$=   {1\over 4 \pi} \   \int_{-\pi}^\pi e^{-x\cosh(t)} \sinh(t)   \int_{-\pi}^{\pi}  \left[ \psi(u)- \psi(-u) \right]  \  {\sin \left((2N+1) (u-t)/2 \right)\over \sin( (u-t) /2)}  du dt.\eqno(3.18)$$
Since $\psi$ is $2\pi$-periodic, we treat  the latter integral with respect to $u$ as follows 

$$  \int_{-\pi}^{\pi}  \left[ \psi(u)- \psi(-u) \right]  \  {\sin \left((2N+1) (u-t)/2 \right)\over \sin( (u-t) /2)}  du $$

$$=  \int_{ t-\pi}^{t+ \pi}  \left[ \psi(u)- \psi(-u) \right]  \  {\sin \left((2N+1) (u-t)/2 \right)\over \sin( (u-t) /2)}  du $$

$$=  \int_{ -\pi}^{\pi}  \left[ \psi(u+t)- \psi(-u-t) \right]  \  {\sin \left((2N+1) u/2 \right)\over \sin( u /2)}  du. $$
Moreover,

$$ {1\over 2\pi} \int_{ -\pi}^{\pi}  \left[ \psi(u+t)- \psi(-u-t) \right]  \  {\sin \left((2N+1) u/2 \right)\over \sin( u /2)}  du - \left[ \psi(t)- \psi(-t) \right] $$

$$=  {1\over 2\pi} \int_{ -\pi}^{\pi}  \left[ \psi(u+t)- \psi(t) + \psi (-t) - \psi(-u-t) \right]  \  {\sin \left((2N+1) u/2 \right)\over \sin( u /2)}  du.$$
When  $u+t > \pi$ or  $u+t < -\pi$ then we interpret  the value  $\psi(u+t)- \psi(t)$ by  formulas

$$\psi(u+t)- \psi(t) = \psi(u+t-2\pi)- \psi(t - 2\pi),$$ 

$$\psi(u+t)- \psi(t) = \psi(u+t+ 2\pi)- \psi(t +2\pi),$$ 
respectively.  Analogously, the value  $\psi(-u-t)- \psi(-t)$  can be treated.   Then   due to the Lipschitz condition (3.15) we have the uniform estimate
for any $t \in [-\pi,\pi]$

$${\left|  \psi(u+t)- \psi(t) + \psi (-t) - \psi(-u-t) \right| \over | \sin( u /2) |}  \le 2C \left| {u\over \sin( u /2)} \right|.$$
Therefore,  owing to the Riemann-Lebesgue lemma

$$\lim_{N\to \infty } {1\over 2\pi} \int_{ -\pi}^{\pi}  \left[ \psi(u+t)- \psi(-u-t)  - \psi(t) + \psi (-t) \right]  \  {\sin \left((2N+1) u/2 \right)\over \sin( u /2)}  du =  0\eqno(3.19)$$
for all $ t\in [-\pi,\pi].$    Besides, returning to (3.18), we estimate the iterated integral 

$$ \int_{-\pi}^\pi e^{-x\cosh(t)} |\sinh(t)  | \int_{ -\pi}^{\pi} \left| \left[ \psi(u+t)- \psi(-u-t)  - \psi(t) + \psi (-t) \right]  \  {\sin \left((2N+1) u/2 \right)\over \sin( u /2)}  \right| du dt$$

$$\le  4 C \int_{0}^\pi e^{-x\cosh(t)}  \sinh(t) dt   \int_{ -\pi}^{\pi}   \left| {u\over \sin( u /2)} \right| du < \infty,\ x >0.$$
Consequently, via  the dominated convergence theorem it is possible to pass to the limit when $N \to \infty$ under the  integral sign, and recalling (3.19), we derive

$$  \lim_{N \to \infty}   {1\over 4 \pi}  \int_{-\pi}^\pi e^{-x\cosh(t)} \sinh(t)  \int_{ -\pi}^{\pi}  \left[ \psi(u+t)- \psi(-u-t)  - \psi(t) + \psi (-t) \right] $$

$$\times  \  {\sin \left((2N+1) u/2 \right)\over \sin( u /2)}  du dt =  {1\over 4 \pi}  \int_{-\pi}^\pi e^{-x\cosh(t)} \sinh(t)  $$

$$ \times \lim_{N \to \infty}  \int_{ -\pi}^{\pi}  \left[ \psi(u+t)- \psi(-u-t)  - \psi(t) + \psi (-t) \right]  \  {\sin \left((2N+1) u/2 \right)\over \sin( u /2)}  du dt = 0.$$
Hence, combining with (3.18),  we obtain  by virtue of  the definition of $\varphi$ and $f$

$$ \lim_{N \to \infty}  S_N(x) =   {1\over 2} \   \int_{-\pi}^\pi e^{-x\cosh(t)}    \left[ \varphi (t)+ \varphi(-t) \right] dt = f(x),$$
where the integral (3.14) converges since $\varphi \in C[-\pi,\pi]$.  Thus we established  (2.36). In order to prove expansion (2.37), we plug the series (3.3) into (2.37)  and recall the representation (2.7) of the kernel $J(y, in, \pi)$ to deduce

$$ {2\over \pi^2} \   \sum_{n=1}^\infty   n  \sinh(\pi n)  K_{in} (x)   \int_0^\infty J (y, in, \pi) f(y) {dy\over y} $$

$$= {2\over \pi^2} \   \sum_{n=1}^\infty   n  \sinh(\pi n)  K_{in} (x)   \int_0^\infty \int_0^\pi e^{-y\cosh(u)} \sinh(u) \sin (n u)   \sum_{m=1}^\infty {a_m \over m}\ K_{im} (y)  du dy. \eqno(3.20)$$
Observing that the interchange of the order of integration and summation is permitted due to the estimates in the proof of Theorem 1, we find the equality

$$ {2\over \pi^2} \int_0^\infty \int_0^\pi e^{-y\cosh(u)} \sinh(u) \sin (n u)   \sum_{m=1}^\infty {a_m\over m}  \ K_{im} (y)  du dy $$

$$= {2\over \pi} \sum_{m=1}^\infty {a_m \over  m \sinh(\pi m)} \int_0^\pi  \sin (n u) \sin (m u) du= {a_n\over n \sinh (\pi n)}.$$
Hence, substituting this  result on the right-hand side of (3.20), we end up with (2.37). 
 
\end{proof} 

{\bf Example 1}.  Let $\psi(u)= \sin(u).$ Then integral (3.14) can be easily calculated by parts to obtain via (2.1) $f(x)= 2 J(x,i,\pi)/x.$  Hence we find the expansion

$$ J(x,i,\pi)  =  {2\over \pi^2} \   \sum_{n=1}^\infty   n  \sinh(\pi n)  J (x, in, \pi)  \int_0^\infty  K_{in} (y) J(y,i,\pi) {dy\over y},\quad x >0.$$

{\bf Example 2}.  Let $\psi $ be the $2\pi$-periodic extension of the function $\psi(u)= u, u \in (-\pi,\pi).$  Then, analogously,  $f(x)= 2 \left(J(x,0,\pi) - \pi e^{-x\cosh(\pi)}\right) /x$  and 

$$ J(x,0,\pi) - \pi e^{-x\cosh(\pi)}  =  {2\over \pi^2} \   \sum_{n=1}^\infty   n  \sinh(\pi n)  J (x, in, \pi) $$

$$\times  \int_0^\infty  K_{in} (y) \left( J(y,0,\pi) - \pi e^{-y\cosh(\pi)} \right) {dy\over y},\quad x >0.$$

Next we give sufficient conditions of the validity of expansions (2.38), (2.39). 

{\bf Theorem 6}. {\it Let $f$ be a complex-valued function on $\mathbb{R}_+$ which is represented by the integral 

$$f(x)= \int_{-\pi}^\pi \sin (x \sinh (u)) \varphi(u) du,\ x >0,\eqno(3.21)$$
where $ \varphi(u) = \psi(u)\cosh(u)$ and $\psi$ is a  $2\pi$-periodic function, satisfying the Lipschitz condition $(3.15)$.
Then expansion $(2.38)$ holds for all $x >0$.  Further, if $f$ can be expanded in terms of the series $(3.3)$
of an arbitrary  sequence $\{a_m\}_{m \in \mathbb{N}}$, satisfying condition $(3.1)$, then the expansion $(2.39)$ holds valid. }

\begin{proof} In fact,  proceeding similar to (3.18) and taking into account formulas (2.8), (2.27),  we have for a partial sum $S_N(x)$ of the series in (2.38)

$$S_{N}(x) =  {1\over \pi} \   \sum_{n=1}^N  \int_{-\pi}^\pi \sin (x \sinh (t)) \cosh(t) \sin (n t) \int_{-\pi}^\pi {\varphi(u)\over \cosh(u)} \sin (n u) dt du$$

$$=   {1\over  4\pi} \   \int_{-\pi}^\pi \sin (x \sinh (t))  \cosh(t)   \int_{-\pi}^{\pi}  \left[ \psi(u)- \psi(-u) \right]  \  {\sin \left((2N+1) (u-t)/2 \right)\over \sin( (u-t) /2)}  du dt.$$
Hence the same ideas and justifications as in the proof of Theorem 5 drive us to the equalities

$$  \lim_{N\to \infty}   S_{N}(x) = {1\over  2} \   \int_{-\pi}^\pi \sin (x \sinh (t))  \left[ \varphi(t) - \varphi(-t) \right] dt  =  \int_{-\pi}^\pi \sin (x \sinh (t)) \varphi(t) dt = f(x),$$
which establish  (2.38).   Next, doing in the manner as in the proof of Theorem 5, we substitute series (3.3) in the integral on the right-hand side of (2.39) and  change the order of integration owing to condition (3.1). Then, appealing to equalities (2.8), (2.27), we find 

$$  {2\over \pi^2} \   \sum_{n=1}^\infty   n  \sinh(\pi n)  K_{in} (x)  \int_0^\infty K_c(y, in, \pi) f(y) {dy\over y} $$

$$= {2\over \pi}  \lim_{N \to \infty}  \sum_{n=1}^N   \sinh \left({\pi n\over 2}\right)  K_{in} (x)   \sum_{m=1}^\infty   {a_m \over   \sinh \left(\pi m/2\right)}\int_0^\pi  \sin (m u)   \sin (n u) du $$

$$=  \sum_{n=1}^\infty a_n  K_{in} (x)  = f(x).$$
Theorem 6 is proved. 

\end{proof}

{\bf Example 3}.  Let $\psi(u)= \sin(u).$ Then integral (3.21) can be easily calculated by parts to obtain via (2.4) $f(x)= 2 \cosh(\pi/2) K_c(x,i,\pi)/x.$  Hence we find the expansion

$$ K_c(x,i,\pi)  =  {2\over \pi^2} \   \sum_{n=1}^\infty   n  \sinh(\pi n)  K_c(x, in,  \pi)  \int_0^\infty  K_{in} (y)  K_c(y,i,\pi) {dy\over y},\quad x >0.$$

Finally we will prove an expansion theorem which concerns equalities (2.40), (2.41). 

{\bf Theorem 7.} {\it Let  $f$ have the following series representation

$$f(\tau)= \tau \sinh\left({\pi \tau\over 2}\right)  \sum_{m=1}^\infty   K_s\left(m, i\tau, \sinh^{-1}(\pi)\right) a_m,\quad \tau >0,\eqno(3.22)$$
where the sequence $\{a_m\}_{m \ge 1}$ satisfies condition $(3.9)$. Then the expansion $(2.40)$ holds for all $x >0$. If, in turn, $f$ is given by series $(3.3)$ under condition $\{a_m\}_{m\ge 1} \in l_1$ then expansion $(2.41)$ takes place for all $x > 0$, where the corresponding integral is understood in the Abel   sense $(2.34)$.}

\begin{proof}  Indeed, recalling  (2.28), (2.29) and (3.10), we plug the series (3.22) into  (2.40) to obtain under condition (3.9) the chain of equalities 

$$ {4\over \pi^2\ } \  x \sinh\left({\pi x\over 2}\right)   \lim_{N\to \infty}  \sum_{n=1}^N  {1\over n} \ K_s\left(n, ix, \sinh^{-1} (\pi) \right) $$

$$\times \int_0^\infty   \cosh\left({\pi \tau\over 2}\right )  K_{i\tau}(n)  f(\tau) d\tau$$

$$= {2\over \pi^2\ } \  x \sinh\left({\pi x\over 2}\right)   \lim_{N\to \infty}  \sum_{n=1}^N  {1\over n} \ K_s\left(n, ix, \sinh^{-1}(\pi) \right) $$

$$\times \int_0^\infty   \tau \sinh\left(\pi \tau\right )  K_{i\tau}(n) \sum_{m=1}^\infty   K_s\left(m, i\tau,  \sinh^{-1}(\pi) \right) a_m d\tau$$

$$= {2\over \pi } \  x \sinh\left({\pi x\over 2}\right) \lim_{N\to \infty}  \sum_{n=1}^N  {1\over n} \ K_s\left(n, ix, \sinh^{-1}(\pi) \right) $$

$$\times \sum_{m=1}^\infty   m  a_m  \int_0^\pi \cos(mu) \cos(nu) du =  x \sinh\left({\pi x\over 2}\right)  \sum_{n=1}^\infty  K_s\left(n, ix, \sinh^{-1}(\pi) \right) a_n = f(x). $$
In order to prove (2.41) we write under conditions of the theorem and equalities (3.10), (3.12)

$${2\over \pi^2 }     \sum_{n=1}^\infty  { K_{ix}(n) \over n}   \int_0^\infty   \tau \sinh\left(\pi \tau\right )  K_s\left(n, i\tau, \sinh^{-1}(\pi)\right)   f(\tau) d\tau$$
 
$$= {4\over \pi^2 }     \sum_{n=1}^\infty  { K_{ix}(n) \over n}   \lim_{\alpha \to {\pi\over 2} -}   \int_0^\infty   \tau \sinh\left({\pi \tau\over 2}\right ) \cosh\left(\alpha \tau \right ) K_s\left(n, i\tau, \sinh^{-1}(\pi) \right) $$

$$\times  \sum_{n=1}^\infty  K_{i\tau}(m) a_m   d\tau = {4\over \pi^2 }     \sum_{n=1}^\infty   K_{ix}(n)  \lim_{\alpha \to {\pi\over 2} -}   \int_0^\infty   \cosh\left(\alpha \tau \right ) \int_0^{\pi} \cos(n u)  \cos \left(\tau \sinh^{-1}(u) \right) du $$

$$\times  \sum_{n=1}^\infty  K_{i\tau}(m) a_m   d\tau = {2\over \pi }     \sum_{n=1}^\infty   K_{ix}(n)  \lim_{\alpha \to {\pi\over 2} -}    \sum_{n=1}^\infty  a_m   \int_0^{\pi}  e^{-m\cos(\alpha) (1+u^2)^{1/2}}  \cos(n u) \cos( m u\sin(\alpha))  du $$

$$=   {2\over \pi }     \sum_{n=1}^\infty   K_{ix}(n)   \sum_{n=1}^\infty  a_m   \int_0^{\pi}    \cos(n u) \cos( m u)  du =  \sum_{n=1}^\infty   K_{ix}(n)  a_n = f(x).$$
where all interchanges of the summation,  integration and the passage to the limit are allowed via the absolute and uniform convergence.   Theorem 7 is proved.

\end{proof}

Let us demonstrate a few examples of the absolutely convergent series and their values which are based on expansions (2.38), (2.40). 

{\bf Example 4}.  Let $f(x)= e^{-x} x^\alpha,\  \alpha > -1.$  Then the integral in (2.38) is calculated via relation (8.4.23.3) in [3], Vol. III, and we find the formula

$$   \sum_{n=1}^\infty   n  \sinh(\pi n) \left| \Gamma( 1+\alpha+ in)\right|^2  K_c(x, in, \pi) =  \pi^{3/2} 2^\alpha  \Gamma\left({3\over 2} +\alpha\right)  x^{\alpha+1} e^{-x},\quad x >0,$$
where $\Gamma(z)$ is the Euler gamma function [3], Vol. III.

{\bf Example 5}.  Let $f(x)= K_\nu(x) x^{\alpha-1},\  \alpha > |\nu|,\ \nu \in \mathbb{R}.$  Then the integral in (2.38) is calculated by formula  (2.16.33.2) in [3], Vol. II,  and we get the equality

$$   \sum_{n=1}^\infty   n  \sinh(\pi n) \left| \Gamma\left( {\alpha+ \nu+ in\over 2}\right) \Gamma\left( {\alpha- \nu+ in\over 2}\right)\right|^2  K_c(x, in, \pi) =   \pi^{2}  2^{2-\alpha}  \Gamma\left(\alpha\right)  x^{\alpha} K_\nu(x),\quad x >0.$$

{\bf Example 6}.  Let $f(x)= \hbox{sech} (\pi x/ 2).$  Then the integral in (2.40) is calculated by formula  (2.16.48.1) in [3], Vol. II,  and we obtain

$$   \sum_{n=1}^\infty   {e^{-n} \over n} \ K_s \left(n, ix, \sinh^{-1} (\pi) \right)=   {\pi\over x \sinh(\pi x) },\quad x >0.$$

{\bf Example 7}.  Let $f(x)=  x \tanh \left(\pi x/ 2\right) \hbox{sech}\left(\pi x/ 2\right).$  Then the integral in (2.40) is calculated by formula  (2.16.48.14) in [3], Vol. II,  and we derive

$$   \sum_{n=1}^\infty   K_0(n) K_s\left(n,  ix, \sinh^{-1} (\pi) \right)=   {\pi^2\over 4} \hbox{sech}^2 \left(\pi x/ 2\right),\quad  x \in \mathbb{R}.$$

\section{Boundary value problem} 

In this section we will employ expansion (2.36) to solve explicitly a Dirichlet boundary value problem in the upper half-plane for the inhomogeneous Helmholtz equation

$$\Delta u - u = h(x,y),\eqno(4.1)$$
where $\Delta$ is the Laplacian,  $(x,y) \in \mathbb{R}_+\times \mathbb{R}_+$, $r= \sqrt{x^2+ y^2}$,

$$h(x,y) =  {2\sinh(\pi)   \over \pi^2 r } e^{- \cosh(\pi) r} \sum_{n=1}^\infty  (-1)^{n+1}  n  \sinh\left( n \cos^{-1} \left({x\over r}\right)  \right)  a_n\eqno(4.2)$$
and $\{a_n\}_{n \ge 1}$ is a suitable sequence.  Writing (4.1) in polar coordinates $(r,\theta),\ r > 0,\ \theta \in [0, \pi]$, it has

$$  {\partial^2u\over \partial r^2}  + {1\over r}  {\partial u\over \partial r} + {1\over r^2}  {\partial^2 u\over \partial \theta^2}-  u$$

$$ = {2 \sinh(\pi)   \over \pi^2 r } e^{- \cosh(\pi) r} \sum_{n=1}^\infty  (-1)^{n+1}  n  \sinh\left( n \theta \right)  a_n.\eqno(4.3)$$

{\bf Theorem 8}. {\it Let  $r >0, \ 0 \le \theta \le \theta_0 < \pi$ and $\{a_n\}_{n \ge 1}$ satisfy the condition

$$\sum_{n=1}^\infty |a_n|\   n^2 e^{\theta_0 n}  < \infty.\eqno(4.4)$$

Then the function

$$u(r,\theta) =  {2\over  \pi^2} \   \sum_{n=1}^\infty   n  \sinh(\theta n)  J (r, in, \pi) a_n\eqno(4.5)$$

is a  solution of equation $(4.3)$, vanishing at infinity.}

\begin{proof}  Indeed, it follows immediately, substituting (4.5) into (4.3) and  fulfilling the differentiation under the series sign due to the absolute and uniform convergence, which is allowed   by virtue of assumption (4.4) and the  estimate  $|J (r, in, \pi) | \le 2 e^{-r}/ n$ in Remark 1.  Then, appealing to the differential equation (2.6) for the incomplete modified Bessel function, we see that $u(r,\theta)$ is a solution of (4.3), which tends to zero  when $r \to \infty$. 

\end{proof}

 {\bf Theorem 9}. {\it Let the sequence $\{a_n\}_{n \ge 1}$ under assumption $(4.4)$ be the Kontorovich-Lebedev transform of discrete variable
 
 $$a_n= \int_0^\infty K_{in} (x) f(x) {dx\over x},\quad n \in \mathbb{N},\eqno(4.6)$$
where $f(x)/x,\ x >0$ satisfies conditions of Theorem $5$. Then  $u(r,\theta)$ given by formula $(4.5)$ solves the inhomogeneous Helmholtz equation $(4.3)$ in the interior of the upper half-plane $(r,\theta),\ r >0, \ 0 < \theta < \pi$, vanishes at infinity  and on the real axis  takes the following prescribed values}

$$  u(r, 0) = 0,  \quad\quad    u(r,\pi) = f(r).\eqno(4.7)$$ 

\begin{proof} The proof follows from  the previous theorem.  To verify boundary conditions (4.7), we appeal to (4.6), Theorem 5 and expansion (2.36). Moreover, $u(r,\theta)$ is continuous on the boundary. In fact, it follows from the uniform convergence with respect to $\theta \in [0,\pi]$ of the series (4.5) by virtue of the Abel test.

\end{proof}

\bigskip
\centerline{{\bf Acknowledgments}}
\bigskip

\noindent The work was partially supported by CMUP, which is financed by national funds through FCT (Portugal)  under the project with reference UIDB/00144/2020.   The author is sincerely indebted to  referees for them  careful reading of the manuscript and constructive comments and suggestions   that greatly improved its form and content.

\bigskip
\centerline{{\bf References}}
\bigskip
\baselineskip=12pt
\medskip
\begin{enumerate}

\item[{\bf 1.}\ ]  D.S. Jones, Incomplete Bessel functions. I,  {\it Proc. Edinb. Math. Soc.} {\bf 50} (2007), N 1,  173-183.

\item[{\bf 2.}\ ] N.N. Lebedev,   Sur un formule d'inversion,  {\it C.R. (Doklady) Acad. Sci. URSS (N.S.)} {\bf 52} (1946),  655-658 (in French).

\item[{\bf 3.}\ ] A.P. Prudnikov, Yu.A. Brychkov and O.I. Marichev, {\it Integrals and Series}. Vol. I: {\it Elementary
Functions}, Vol. II: {\it Special Functions}, Gordon and Breach, New York and London, 1986, Vol. III : {\it More special functions},  Gordon and Breach, New York and London,  1990.

\item[{\bf 4.}\ ] S. Yakubovich, {\it Index Transforms}, World Scientific Publishing Company, Singapore, New Jersey, London and
Hong Kong, 1996.

\end{enumerate}

\vspace{5mm}

\noindent S.Yakubovich\\
Department of  Mathematics,\\
Faculty of Sciences,\\
University of Porto,\\
Campo Alegre st., 687\\
4169-007 Porto\\
Portugal\\
E-Mail: syakubov@fc.up.pt\\

\end{document}